\documentclass{amsart}
\usepackage{amsmath,amssymb, amsthm, color,hyperref}
\usepackage{tikz}
\usetikzlibrary{calc,through,intersections, arrows, decorations.markings, matrix}
\tikzset{>=latex}

\theoremstyle{plain}
\newtheorem{theorem}{Theorem}[section]

\theoremstyle{definition}

\newcommand{\integers}{\ensuremath{\mathbb Z}} 
\newcommand{\CC}{\mathbf{C}}
\newcommand{\TT}{\mathbf{T}}

\newcommand{\DD}{\mathbf{D}}

\renewcommand{\top}{\mathbf{Top}}

\renewcommand{\hom}{\text{Hom}}

\newcommand{\sSet}{\mathbf{sSet}}
\newcommand{\alg}{\mathbf{Alg}}

\newcommand{\map}{\text{Map}}

\newcommand{\hocolim}{\mathrm{hocolim}}

\begin{document}

\title{Spaces with Complexity One}
\author{Alyson Bittner}
\begin{abstract}
An $A$-cellular space is a space built from $A$ and its suspensions,  analogously to the way that $CW$-complexes are built from $S^0$ and its suspensions. The  $A$-cellular approximation of a space $X$
is an $A$-cellular space $CW_{A}X$ which is  closest to $X$ among all $A$-cellular spaces. 
The $A$-complexity of a space $X$ is an ordinal number that quantifies how difficult it is to build an 
$A$-cellular approximation of  $X$.  In this paper, we study spaces with low complexity. In particular we show that if 
$A$ is a sphere localized at a set of primes then the $A$-complexity of each space $X$ is at most 1. 

\end{abstract}

\maketitle
\section{Introduction}

Let $A$ be a pointed CW-complex. An $A$-cellular space is a space built out of copies of $A$ via iterated pointed homotopy colimit construction. Given a pointed space $X$, the $A$-cellular approximation of $X$ is an $A$-cellular space 
$CW_{A}X$ equipped with a map $CW_{A}X\to X$,  such that the induced map of pointed mapping complexes 
$\map_{\ast}(A, CW_{A}X) \to \map_{\ast}(A, X)$ is a weak equivalence.  For example, $S^{0}$-cellular spaces are CW-complexes and an $S^{0}$-cellular approximation of a space X is the CW-approximation of X. 
Farjoun \cite{Farjoun1995} showed that the $A$-cellular approximation of $X$ exists for any $A$ and $X$ and that it is unique up to homotopy equivalence. 
The space $CW_{A}X$ is the best possible approximation of $X$ in the class of $A$-cellular spaces:
for any $A$-cellular space $Y$ and a map $Y\to X$ there exists a map $Y \to CW_A X$, unique up to homotopy, 
such that the following diagram is homotopy commutative.
\begin{center}
$$
\begin{tikzpicture}
\matrix (m) 
[matrix of math nodes, row sep=2em, column sep=3em, text height=1.5ex, text depth=0.25ex]
{
Y &  CW_{A}X\\
 &   X\\
};
\path[->]
(m-1-2)  edge  (m-2-2)
(m-1-1)   edge (m-2-2)
;
\path[->, dashed]
(m-1-1) edge  (m-1-2)
;
\end{tikzpicture}
$$ 
\end{center}

In \cite{ChacholskiDwyerIntermont2002} Chach\'olski, Dwyer, and Intermont introduced the concept of the $A$-complexity 
of a space $X$, which is the minimum ordinal number of homotopy colimits necessary to produce $CW_{A}(X)$ from copies of $A$. More precisely, starting with the full subcategory  $\CC_{0}$ of pointed spaces whose  objects have  the homotopy type of  a retract of wedges of $A$, one can construct an increasing  chain of categories  indexed by ordinal numbers $\alpha$, where  $\CC_{\alpha}$ is the category of all spaces of the homotopy type of $\hocolim_{\DD}F$ for some small 
category $\DD$ and a functor $F:\DD\to \bigcup_{\beta < \alpha} \CC_{\beta}$ and their retracts. The  the $A$-complexity of $X$ is the minimal ordinal $\kappa_{A}(X)$ such that $X\in \CC_{\kappa_{A}(X)}$.

In the case of $A = S^{0}$ the the category $\CC_{1}$ consists of all space of a homotopy type of a CW-complex, which 
implies that $\kappa_{S^{0}}(X) = 1$ for any space  $X$. This is however not the case for other choices of $A$, and in 
general the complication of a space need not even be finite.  For example, if $A = M(\integers/p,n)$  is a Moore space for some $n\ge 1$  and  a prime $p$  and  $X = M(\integers/p^{\infty}, n+1)$  then $\kappa_{A}(X) = \omega$ (see  \cite[Proposition 9.3]{ChacholskiDwyerIntermont2002}). 
By contrast,  Chach\'olski, Dwyer and Intermont observed that if $A=S^{n}$ then $\kappa_{A}(X) \leq 2$ for all $X$  and suggested that it should be possible to lower this bound 
and  show that $\kappa_{A}(X) \leq 1$  \cite[\S 9.3]{ChacholskiDwyerIntermont2002}. The goal of this note is to show 
that this is indeed the case, and that this result  holds in even greater generality:   

\begin{theorem}
\label{MAIN THM}
If $A$ is a sphere, or a sphere localized at a set of primes then $\kappa_{A}(X) \leq 1$ for all spaces $X$.
\end{theorem}

\section*{Acknowledgements}

I want to express my appreciation to my advisor, Bernard Badzioch, for his helpful comments which led to the final version of paper. Additionally, I would like to thank Matthew Sartwell for several useful conversations.

\section{Algebras of mapping spaces}

Our proof of Theorem \ref{MAIN THM} will be based on the following algebraic description of mapping space 
obtained by Sartwell  \cite{SartwellThesis}. For  spaces $X, Y$ let $\map_{\ast}^{\Delta}(X, Y)$ denote the 
simplicial mapping complex of pointed maps $X\to Y$. Given a pointed CW-complex $A$ let $\TT^{A}$ denote 
a simplicial category 
on objects $T_{n}$ for $n=0, 1, \dots$, and such that 
$$\hom_{\TT^{A}}(T_n, T_m)=\map_*^{\Delta}\left(\bigvee^m A, \bigvee^n A\right)$$
Notice that in this category $T_{n}$ is the $n$-fold product of $T_{1}$. A $\TT^{A}$-algebra is 
a product preserving simplicial functor $\TT^{A}\to \sSet_*$ where $\sSet_*$ is the category of pointed
simplicial sets.  Let $\alg^{\TT^{A}}$ denote the category of all $\TT^{A}$-algebras with natural transformations as morphisms.

Notice that any pointed space $X$ defines a $\TT^{A}$-algebra 
$\Omega^{A}(X)$ such that  
$$\Omega^{A}(X)(T_{n}) = \map_{\ast}^{\Delta}\left(\bigvee^{n} A, X\right)$$ 
The resulting functor $\Omega^{A}\colon \top_{\ast} \to \alg^{\TT^{A}}$ has a left adjoint $B^{A}$.

The category $\alg^{\TT^{A}}$ can be equipped with a model category structure where weak equivalences and 
fibrations are defined as objectwise weak homotopy equivalence and  Serre fibrations respectively.  Denote also by 
${\mathbf R}^{A}\top_{\ast}$ the category $\top_{\ast}$ taken with the model category structure where 
fibrations are Serre fibrations and a weak equivalences are maps $f\colon X \to Y$ that induce a weak equivalences
$f_{\ast}\colon \map_{\ast}^{\Delta}(A, X) \to \map_{\ast}^{\Delta}(A, Y)$.  In other words,  ${\mathbf R}^{A}\top_{\ast}$
is obtained by taking the right Bousfield localization of the usual model category structure on $\top_{\ast}$ with respect 
to the space $A$ as in \cite[5.1.1]{Hirschhorn2009}.  Sartwell showed that the following holds:

\begin{theorem}\cite{SartwellThesis}\label{Adjunction}
Let $A$ be a sphere or a sphere localized at a set of primes. The adjoint pair 
$$B^{A}\colon \alg^{\TT^{A}} \rightleftarrows {\mathbf R}^{A}\top_{\ast} \colon \Omega^{A}$$
is a Quillen equivalence. 
\end{theorem}

\section{Proof of Theorem \ref{MAIN THM}}

Directly from the definition of an $A$-cellular approximation it follows that if spaces $X$ and $Y$ are weakly equivalent in $R^A\top$ then $CW_AX$ and $CW_AY$ are weakly homotopy equivalent, and so $\kappa_A(X)=\kappa_A(Y)$. Let $A$ be a space as in Theorem \ref{MAIN THM}.  By Theorem \ref{Adjunction} for any space $X$ we get a weak equivalence 
$X \simeq B^AQ\Omega^AX$  in $R^A\top$, where $Q\Omega^AX$ denoted a cofibrant replacement in 
the category $\alg^{\TT^{A}}$.  Therefore it is only necessary to show that for any space $X$,
$\kappa_A( B^AQ\Omega^AX) \leq 1$.

The algebra $Q\Omega^AX$ can be described more explicitly as follows. Let  $U:\alg^{\TT^A}\to\sSet_*$ denote 
the forgetful functor, $U(\Phi)=\Phi(T_1)$.  This functor has a left adjoint $F:\sSet_*\to\alg^{\TT^A}$. 
For $\Phi \in \TT^A$-algebra, let  $FU_{\bullet}\Phi$ be the simplicial $\TT^A$-algebra defined by adjoint pair 
$(F, U)$. Explicitly:
$$FU_n\Phi:=(FU)^{n+1}\Phi$$
By \cite[2.3.3]{SartwellThesis} the natural map $|FU_{\bullet}\Phi| \to \Phi$ is a cofibrant replacement of $\Phi$ in 
the category $\alg^{\TT^{A}}$. 

In view of this fact we need to show that for any space $X$ the $A$-complexity of  $B^A|FU_{\bullet}\Omega^AX|$ is at most one. The left adjoint $B^A$ can be described explicitly as the coend $B^A(\Phi)=\int^{T_n\in T^A}\bigvee_nA\wedge |\Phi(T_n)|$, which is in a similar manner as done in \cite[3.4]{BadziochChungVoronov2007} but then composing with objectwise singularization and realization to obtain the desired adjunction \cite[18.3.10]{Hirschhorn2009}. 
The functor $B^AF\colon \sSet_* \to \top_{\ast}  $ is left adjoint to $U\Omega^A$, and since 
$U\Omega^A(X)=\map^{\Delta}_*(A, X)$, thus  $B^AF(Y)= A\wedge |Y|$. 
This gives:
\begin{align*}
B^A|FU_{\bullet}\Omega^AX|  \simeq\  |B^AFU_{\bullet}\Omega^AX|\simeq& \hocolim_{n\in \Delta^{op}} A\wedge |(FU)^{n}\Omega^{A}X(T_{1})| \\
 \simeq\  & \hocolim_{(n, m)\in \Delta^{op}\times \Delta^{op}} A\wedge ((FU)^{n}\Omega^{A}X(T_{1}))_{m} \\
 \simeq\  & \hocolim_{(n, m)\in \Delta^{op}\times \Delta^{op}}  \bigvee^{\sigma\in((FU)^{n}\Omega^{A}X(T_{1}))_{m}} A \\
\end{align*}

\bibliographystyle{plain}
\bibliography{MyPaperResources}

\end{document}